\documentclass[final, 3p, times]{elsarticle}

\usepackage{hyperref, float, amsmath, bm, listings, color, footnote, caption, xparse}
\usepackage{amsfonts, amssymb, comment, graphicx, subcaption, listings, tikz}
\usetikzlibrary{shapes, arrows}
\usepackage[linesnumbered,ruled,vlined]{algorithm2e}
\usepackage[bottom]{footmisc}
\tikzstyle{startstop} = [rectangle, rounded corners, minimum width=3cm, minimum height=1cm, text centered, draw=black, fill=red!30]
\tikzstyle{io} = [trapezium, trapezium stretches=true, trapezium left angle=70, trapezium right angle=110, minimum width=3cm, minimum height=1cm, text centered, draw=black, fill=blue!30]
\tikzstyle{process} = [rectangle, minimum width=1cm, minimum height=1cm, text centered, text width=4cm, draw, fill=orange!30]
\tikzstyle{decision} = [diamond, minimum width=1cm, minimum height=1cm, text centered, draw=black, fill=green!30]
\tikzstyle{arrow} = [thick,->,>=stealth]

\newcommand{\lexint}{\texttt{LeXInt}\xspace}

\hypersetup{colorlinks,
            linkcolor = {red!50!black},
            citecolor = {red!10!black},
            urlcolor  = {green!10!black}}

\definecolor{codegreen}{rgb}{0.1, 0.1, 0.9}
\definecolor{codegray}{rgb}{0, 0.5, 0}
\definecolor{codepurple}{rgb}{0.6, 0.1, 0.9}
\definecolor{codeorange}{rgb}{1.0, 0.5, 0.1}
\definecolor{backcolour}{rgb}{0.99, 0.99, 0.99}

\lstdefinestyle{CodeColour}{
	backgroundcolor=\color{backcolour},   
	commentstyle=\color{codegreen},
	keywordstyle=\color{codeorange},
	numberstyle=\tiny\color{codegray},
	stringstyle=\color{codepurple},
	basicstyle=\ttfamily\footnotesize,
	breakatwhitespace=false,         
	breaklines=true,                 
	captionpos=b,                    
	keepspaces=true,                 
	numbers=left,                    
	numbersep=5pt,                  
	showspaces=false,                
	showstringspaces=false,
	showtabs=false,                  
	tabsize=4
}

\journal{SoftwareX}

\begin{document}

\begin{frontmatter}



\title{\lexint: GPU-accelerated Exponential Integrators package}

\author{Pranab J. Deka\footnote{Now at KU Leuven, Belgium}} \ead{pranab.deka@kuleuven.be} 
\author{Alexander Moriggl}
\author{Lukas Einkemmer\corref{lod1}} \ead{lukas.einkemmer@uibk.ac.at} 

\cortext[lod1]{Corresponding author}

\affiliation{organization = {Department of Mathematics, University of Innsbruck},
            city = {Innsbruck},
            postcode = {6020}, 
            country = {Austria}}

\begin{abstract}
We present an open-source \texttt{CUDA}-based package, for the temporal integation of differential equations,  that consists of a compilation of exponential integrators where the action of the matrix exponential or the $\varphi_l$ functions on a vector is approximated using the method of polynomial interpolation at Leja points. Using a couple of test examples on an NVIDIA A100 GPU, we show that one can achieve significant speedups using \texttt{CUDA} over the corresponding CPU code. \lexint, written in a modular format, facilitates integration into existing software packages (written in \texttt{C++} or \texttt{CUDA}), and can be used for temporal integration of differential equations.

\end{abstract}


% \begin{highlights}
% \end{highlights}


\end{frontmatter}


%%%%%%%%%%%%%%%%%%%%%%%%%%%%%%%%%%%%%%%%%%%%%%%%%%%%%%%%%%%%%%%%%%%%%%%%%%%%%%%%%%

\begin{table*}[h]
\captionof{table}{Code metadata}
\small
\hspace*{-0.75cm}
\centering
\begin{tabular}{ll}
    \hline
    Code Version   & 1.0.1  \\
    Permanent link to code/repository & \url{https://github.com/Pranab-JD/LeXInt}  \\
    Legal Code License &  MIT \\
    Code versioning system used & Git \\
    Software code languages, tools, and services used & C++, CUDA \\
    Compilation requirements, operating environments \& dependencies & NVIDIA GPU, CUDA11.2+, g++ and nvcc compilers \\
    Support email for questions & \href{mailto:pranab.deka@kuleuven.be}{pranab.deka@kuleuven.be} \\
    \hline
\end{tabular}
\end{table*}


\section{Motivation and significance}

Exponential integrators have emerged as viable alternatives to implicit or semi-implicit integrators in numerically treating ordinary and partial differential equations (ODEs/PDEs). Owing to their ability to choose significantly larger time step sizes and their superior stability properties as compared to implicit or semi-implicit ones, they are well-suited for stiff problems, i.e., problems that involve physical processes with both slow- and fast-changing dynamics, resulting in widely varying time scales.  Furthermore, they have the added advantage that they can solve the linear (stiff) part of an equation \textit{exactly} (in time). Several classes of exponential integrators have been developed over the years - Lawson integrators \cite{Lawson76}, exponential time differencing \cite{Beylkin98, Cox02}, integrating factor \cite{Boyd89, Fornberg99}, exponential Rosenbrock (EXPRB, \cite{Caliari09, Ostermann10}), and exponential propagation iterative methods of Runge--Kutta type (EPIRK, \cite{Tokman06, Tokman12}). Let us consider the abstract PDE,
\begin{equation}
    \frac{\partial u}{\partial t} = f(u) = \mathcal{A} \, u + g(u),
    \label{eq:pde}
\end{equation}
where $\mathcal{A}\,u$ is the linear (usually stiff) term and $g(u)$ is the nonlinear (usually nonstiff) term of some nonlinear function $f(u)$. Exponential integrators require approximating the exponential of the matrix ($\exp(\mathcal{A})$) or exponential-like $\varphi_l(\mathcal{A})$ functions. These $\varphi_l(\mathcal{A})$ functions are defined recursively as
\[\varphi_{l+1}(\mathcal{A}) = \frac{1}{\mathcal{A}} \left(\varphi_l(\mathcal{A}) - \frac{1}{l!} \right), \; \, l\geq 0, \; \, \text{with} \; \, \varphi_0(\mathcal{A}) = \exp(\mathcal{A}).\]

The simplest EXPRB method, the second-order exponential Rosenbrock--Euler integrator, was proposed by \citep{Pope63}. Exponential integrators did not receive a lot of attention from the scientific community in the early 60s owing to the hefty computational expenses of approximating the matrix exponential for large problems \citep{Tokman06, Tokman10, Ostermann10}. The advent of Krylov-subspace methods \cite{Vorst87}, based on the Arnoldi algorithm \cite{Arnoldi1951}, enabled the efficient computation of the matrix exponential applied to a vector for large problems. This has managed to successfully reignite the interest of the scientific community in exponential integrators. Krylov-based exponential integrators have been shown to be highly competitive or even surpass the state of the art implicit integrators, especially for stiff problems \cite{Hochbruck09, Ostermann10, Tokman06, Tokman13, Tokman16, Einkemmer17}. Some of the publicly available Krylov-based exponential integrators software include \texttt{EXPOKIT} \cite{Sidje98}, \texttt{phipm} \cite{Niesen09}, \texttt{phipm\_simul\_iom} \cite{Luan19}, and \texttt{EPIC} \footnote{\url{https://faculty.ucmerced.edu/mtokman/\#software}}.

Approximating the $\varphi_l$ functions applied to a vector as polynomials on the set of Leja points has been proposed as an alternative to the Krylov-based methods \cite{Caliari04, Bergamaschi06, Caliari07a, Caliari07b} to compute the matrix exponential and the related $\varphi_l$ functions. The main advantages of the Leja-based approach over Krylov are as follows:
\begin{itemize}

    \item Polynomial interpolation at Leja points, in many scenarios, can outperform Krylov-based methods \citep{Bergamaschi06, Caliari07a, Deka22b} owing to the simplicity of the algorithm as well as its implementation.

    \item Whilst Krylov-based methods may require several basis vectors to be stored in memory (which can be a serious impediment for large-scale problems), the Leja method requires the storage of only one input and one output vector (for computing the matrix exponential and for single-stage exponential integrators). This is also favourable for modern computer architectures, such as graphical processing units (GPUs), where memory tends to be limited.

    \item The need for the evaluation of inner products in Krylov-based methods may result in performance loss when scaling such methods to large supercomputers or massively parallel systems.
    
\end{itemize}

With the release of the \texttt{Python}-based open-source software, \textbf{Le}ja interpolation for e\textbf{X}ponential \textbf{Int}egrators or \lexint \cite[\url{https://github.com/Pranab-JD/LeXInt}]{Deka22_lexint}, we provided a user-friendly framework for exponential integrators and highly-efficient novel algorithms that are used in such solvers. We aimed for the scientific community to get started exploring such novel methods in a high-level language. 

However, it is well-known that large-scale simulations of complex physical processes such as galaxy evolution, ocean dynamics, or plasma processes require the use of low-level programming languages owing to their superior processing speeds. Parallelising large-scale problems over multiple processors is desirable due to the hefty memory requirements and to achieve significant improvement in speed in the simulations. The ever increasing need and desire for high-resolution long-time simulations have popularised the use of accelerators, such as GPUs, to achieve even larger boosts. High memory bandwidths make GPUs highly attractive for memory-bound scientific problems. Furthermore, on average, GPUs are less energy-demanding than CPUs, thereby making them more desirable. Compute Unified Device Architecture (\texttt{CUDA}), developed by NVIDIA, is a software that enables direct access of the high throughput of the GPUs. \texttt{CUDA} \footnote{\url{https://docs.nvidia.com/cuda/doc/index.html}} provides excellent documentation that enables easy usage and includes tools for compiling, profiling, and debugging. 

In this work, we present a \texttt{CUDA} (and \texttt{C++}) implementation of \lexint, which, to the best of our knowledge, is the first open-source \texttt{CUDA}-based software for exponential integrators. It fully exploits the high throughput offered by GPUs by performing over $90\%$ of the computation on the GPUs and keeping the number of transfers between the host and the device memory to a bare minimum. This \texttt{C++/CUDA} version of the code has been well-optimised for memory-bound problems that often result from the discretisation of PDEs. \lexint consists of a wide range of exponential integrators, namely Rosenbrock--Euler \cite{Pope63}, EXPRB32 \cite{Caliari09, Ostermann10}, EXPRB43 \cite{Caliari09, Ostermann10}, EXPRB42 \cite{Luan17}, EXPRB53s3 \cite{Luan14}, EXPRB54s4 \cite{Luan14}, EPIRK4s3 \cite{Tokman17a, Tokman17b}, EPIRK4s3A \cite{Tokman16}, EPIRK4s3B \cite{Tokman16}, and EPIRK5P1 \cite{Tokman12}. The $\varphi_l$ functions, in each of these exponential integrators have been approximated using the Leja interpolation method. Higher-order exponential integrators typically require the computation of one or more internal stages, which in turn requires the storage of a few additional vectors. The number of auxiliary vectors (temporary data) needed to be stored in memory has been optimised to a bare minimum for each of the integrators. \lexint has been developed to be used as a black box throughout the simulations. The simulation data need to be transferred over from the device to the host memory only at the very end of the simulations. However, the user is completely free to have the output transferred over to the host memory and written to files at intermediate stages should they choose to do so. Of course, this will induce an overhead in the overall computational cost. Let us clearly state that the current version of \lexint does not support MPI parallelisation, and hence cannot run on multiple GPUs. We acknowledge that this is crucial for running large-scale simulations over multiple nodes. Facilitating MPI-parallelisation to use multiple GPUs would be a subject of future work.

This article is organised as follows: in Sec. \ref{sec:software}, we describe the memory allocation and provide details as to how to invoke the exponential integrators and the Leja interpolation method. In Sec. \ref{sec:eg}, we illustrate the performance of the \texttt{CUDA}-based \lexint and contrast its performance with \texttt{C++}. Finally, in Sec. \ref{sec:impact}, we highlight the potential impact of this work. 

%%%%%%%%%%%%%%%%%%%%%%%%%%%%%%%%%%%%%%%%%%%%%%%%%%%%%%%%%%%%%%%%%%%%%%%%%%%%%%%%%%

\section{Software description}
\label{sec:software}

\lexint consists of a collection of exponential integrators along with a range of helper functions, such as computing the spectral radius using the power iteration method \cite{Deka22b, Deka22_lexint}. To facilitate the usage of these functions, we define the \texttt{struct Leja}, in \texttt{Leja.hpp}, which collects these functions and handles the memory management for each of the exponential integrators and computation of the action of the Jacobian on the relevant vector. \lexint is a header-only library, i.e., no installation is required. The line \lstinline{#include "./LeXInt/CUDA/Leja.hpp"}, with the proper path to \lexint (depending on the specific directory where the user downloads the package), should be added to the appropriate file within their code. A flowchart illustrating the code structure is shown in Fig. \ref{flow:lexint_structure}.
\begin{figure}[t]
    \centering
    \begin{tikzpicture}[node distance=1.75cm]
        
        \node[draw, align=center] (start) [startstop] {Invoke \lexint \\ (\texttt{\#include ./LeXInt/CUDA/Leja.hpp})};
        \node (pro1) [process, below of=start] {Approximate most dominant eigenvalue ($c, \gamma$)};
        \node (dec1) [decision, below of=pro1, yshift = -1cm] {Linear DE};
        \node (dec2) [decision, right of=dec1, xshift = 7.5cm] {Nonlinear DE};
        \node (dec3) [decision, below of=dec1, yshift = -1.5cm] {Homogenous};
        \node (dec4) [decision, right of=dec3, xshift = 3cm] {Nonhomogenous};
        \node (pro2) [process, below of=dec3, yshift = -1.5cm] {\texttt{real\_Leja\_exp()}};
        \node (pro3) [process, below of=dec4, yshift = -1.5cm] {\texttt{real\_Leja\_phi\_nl() }};
        \node[draw, align=center] (pro4) [process, below of=dec2, yshift = -4.75cm] {\textit{choose desired \\ exponential integrator}};
        
        \draw [arrow] (start) -- (pro1);
        \draw [arrow] (pro1) -- (dec1);
        \draw [arrow] (pro1) -| (dec2);
        \draw [arrow] (dec1) -- (dec3);
        \draw [arrow] (dec1) -| (dec4);
        \draw [arrow] (dec3) -- (pro2);
        \draw [arrow] (dec4) -- (pro3);
        \draw [arrow] (dec2) -- (pro4);

    \end{tikzpicture}
    \caption{A flowchart demonstrating the structure and usage of \lexint.}
    \label{flow:lexint_structure}
\end{figure}

Integrating Eq. \eqref{eq:pde} using \lexint requires the user to define a function that computes the right-hand side (rhs), i.e. $f(u)$. This function is expected to take exactly two arguments - the first parameter is the pointer to the input data, whilst the second is the pointer to the output data, where the data for both the input and output have to lie contiguous in memory. If the user-specified rhs function consists of additional parameters, one could construct a \texttt{class} and have these supplementary parameters as member variables of this \texttt{class} (e.g., Listing \ref{alg:RHS}).

%The \texttt{RHS} operator (of type \texttt{void}) is expected to defined in a way such it takes one-dimensional vector of size $N$ and yields an one-dimensional vector of the same size. If the RHS function consists of additional parameters, one could potentially construct a \texttt{struct} and have these parameters as localised to this \texttt{struct} (see Algorithm \ref{alg:RHS}). They may also be defined as arguments of a parent class. In the test examples provided in \lexint, the parameters such as the number of grid points, advection velocity (wherever applicable), and grid resolution are defined in the parent class \texttt{Problems\_2D}. In the case of diffusion in the presence of a source, the location of the source is defined in the function \texttt{Dif\_Source\_2D} which is a part of the \texttt{struct} \texttt{RHS\_Dif\_Source\_2D} derived from the parent class \texttt{Problems\_2D}.
\begin{lstlisting}[language = C++, caption = Functor to add parameters to the rhs of the equation., label = alg:RHS]
struct rhs
{
    rhs(*args)

    void operator(input, output)
    {
        RHS(input, output, *args)
    }
}
\end{lstlisting}

%%% --------------------------------------------------------- %%%

\subsection{Polynomial Interpolation at Leja points}

Let us assume that $\mathbb{K}$ is a compact set and $\mathbb{K} \subset \mathbb{C}$, where $\mathbb{C}$ is the complex plane. A set of Leja points \cite{Edrei1940, Leja1957} can be defined recursively as \[\prod_{k = 0}^{j - 1} |z_j - z_k| = \text{max}_{z \in \mathbb{K}}\prod_{k = 0}^{j - 1} |z - z_k|,\] where $z \in \mathbb{K}$ and $j = 1, 2, 3 \hdots$. Conventionally, $|z_0|$ is chosen as $\text{max}_{z \in \mathbb{K}} |z|$. We propose the use of Leja points, for computing the matrix exponential and the $\varphi$ functions, owing to its recursive nature. The well-known Chebyshev points, for example, do not exhibit this property - the interpolation of a polynomial at $m + 1$ Chebyshev points necessitates the re-computation of the $\varphi_{l} $functions at the previously computed $m$ nodes. However, since Leja points are recursive in nature, using $m + 1$ Leja points needs only one extra computation, and the computation at the previous $m$ nodes can be reused. This results in significant computational savings.

The method of polynomial interpolation at Leja points \cite{Leja1957, Baglama98, Reichel90} to compute the matrix exponential and the $\varphi_l(z)$ functions was proposed in \cite{Caliari04, Bergamaschi06}, and has subsequently been shown to be highly competitive with the traditionally used Krylov-based methods \cite{Caliari04, Deka22b, Deka23b}. This method has been described theoretically in great detail in \cite{Caliari04, Caliari14}, and we have described the working algorithm in our previous work \cite{Deka22a, Deka22b, Deka22_lexint, Deka23a, Deka23b}. Here, we present a gist of the method: Let us suppose that we want to compute the exponential of a matrix $\mathcal{A}$ times a vector $y_0$, i.e., $\exp(\mathcal{A})y_0$. Let us suppose that the eigenvalues of $\mathcal{A}$ satisfy $\alpha \leq \mathrm{Re} \; \sigma(\mathcal{A}) \leq \nu \leq 0 $ and $-\beta \leq \mathrm{Im} \; \sigma (\mathcal{A}) \leq \beta$, where $\sigma(\mathcal{A})$ is the spectrum of $\mathcal{A}$; $\alpha$ and $\nu$ are the smallest and largest real eigenvalues respectively, and $\beta$ is the largest, in modulus, imaginary eigenvalue. We construct an ellipse, with semi-major axis $a$ and semi-minor axis $b$, encompassing all eigenvalues of $\mathcal{A}$.  Let $c$ and $\gamma$ be the midpoint point and the distance between the two foci of the ellipse, respectively. To compute the matrix exponential, we interpolate the function $\exp(c + \gamma \xi)$ on a set of precomputed Leja points ($\xi$) in the interval $[-2, 2]$. The ${(m+1)}^\mathrm{th}$ term of the polynomial $p$ is given by
\begin{align}
    p_{m+1} & = p_m + d_{m+1} \, y_{m+1}, \nonumber \\
    y_{m+1} & = y_m \times \left(\frac{z - c}{\gamma} - \xi_{m} \right).
    \label{eq:leja}
\end{align}
The polynomial is initialised as $p_0 = d_0 \, y_0$, where $y_0$ is the vector that is applied to the matrix exponential. To approximate the $\varphi_k$ function, one needs to replace $\exp(c + \gamma \xi)$ with $\varphi_k(c + \gamma \xi)$.

%%% --------------------------------------------------------- %%%

\subsection{Code structure (``external")}

In this subsection, we provide a high-level description of the code that would enable a user to use \lexint as a temporal integration package for their software.

\subsubsection*{Leja interpolation}
In \lexint, the functions \texttt{real\_Leja\_exp}, \texttt{real\_Leja\_phi\_nl}, and \texttt{real\_Leja\_phi} carry out the polynomial interpolation at Leja points. Listings \ref{alg:leja_exp_ext} and \ref{alg:leja_phi_nl_ext} correspond to the ``external" structure of the functions \texttt{real\_Leja\_exp} and \texttt{real\_Leja\_phi\_nl}, respectively, i.e, the user needs to provide the arguments only to these functions. Here, \texttt{RHS} is the user-defined function that computes the rhs or $f(u)$ for a given problem. \texttt{u\_input} and \texttt{u\_output} are the input and the output vectors, respectively. It should be noted that data containers in multidimensional problems are expected to lie contiguous in memory, else the performance improvements on GPUs will be severely limited. \texttt{c} and \texttt{Gamma} correspond to $c$ and $\gamma$ in Eq. \ref{eq:leja}, respectively. \texttt{rtol} and \texttt{atol} are the user-defined relative and absolute tolerances, respectively, (converging condition: \texttt{rtol} $\times$ $||$\texttt{polynomial}$||$ + \texttt{atol}), \texttt{dt} is the time step size, \texttt{iters}, which is an output parameter, corresponds to the number of Leja iterations for a particular time step, and \texttt{GPU} = 1 if GPU support is enabled, else, \texttt{GPU} = 0. When GPU support is enabled, \texttt{cublas\_handle} is used to invoke \texttt{cublas}, which is used to compute the l2 norm  of the relevant vectors (l2 norm normalised to $\sqrt{N}$ is used for error estimation). This handle is defined in \texttt{Leja.hpp}.

\begin{minipage}{0.95\linewidth}
\begin{lstlisting}[language = C++, caption = Leja interpolation of the matrix exponential applied to vector `u' (``external" function). The user needs to provide the arguments to this function., label = alg:leja_exp_ext]
//? Real Leja Exp (Homogenous linear equations - matrix exponential)
template <typename rhs>
void real_Leja_exp(rhs& RHS,                   //? RHS function
                   double* u_input,            //? Input state variable(s)
                   double* u_output,           //? Output matrix exponential times 'u'
                   double c,                   //? Shifting factor
                   double Gamma,               //? Scaling factor
                   double rtol,                //? Relative tolerance (normalised)
                   double atol,                //? Absolute tolerance
                   double dt,                  //? Step size
                   int& iters,                 //? # of iterations needed to converge
                   bool GPU,                   //? false (0) --> CPU; true (1) --> GPU
                   )
\end{lstlisting}

\begin{lstlisting}[language = C++, caption = Leja interpolation of the $\varphi$ function applied to \texttt{interp\_vector} (``internal" function).. The user needs to provide the arguments to this function., label = alg:leja_phi_nl_ext]
//? Real Leja Phi NL (Nonhomogenous linear equations)
template <typename rhs>
void real_Leja_phi_nl(
            rhs& RHS,                         //? RHS function
            double* u_input,                  //? Input state variable(s)
            double* u_output,                 //? Output applied to phi function
            double (* phi_function) (double), //? Phi function
            double c,                         //? Shifting factor
            double Gamma,                     //? Scaling factor
            double rtol,                      //? Relative tolerance (normalised)
            double atol,                      //? Absolute tolerance
            double dt,                        //? Step size
            int& iters,                       //? # of iterations needed to converge
            bool GPU                          //? false (0) --> CPU; true (1) --> GPU
            ) 
\end{lstlisting}
\end{minipage}

The user has to provide the RHS function (in the desired template) and the pointers to the input (\texttt{u\_input}) and output (\texttt{u\_output}) vectors. As the interpolation functions are of type \texttt{void}, the user has to assign the required amount of memory for the output vector. This can be done as follows:

\begin{minipage}{0.95\linewidth}
\begin{lstlisting}[language = C++, caption = Define and assign memory to the output array.]
// Pointer to the user-desired name for the array
double* u_output;

// Size of the array
size_t N_size = N * sizeof(double);

// Assign memory on C++ code
u_output = (double*)malloc(N_size);

// Assign memory on CUDA code
cudaMalloc(&u_output, N_size);
\end{lstlisting}
\end{minipage}

The scaling and shifting factors are to be computed before starting the time loop. The desired time step size and tolerance are to be provided by the user. There is the choice of running the program on CPU (\texttt{bool GPU = 0}) or NVIDIA GPU  (\texttt{bool GPU = 1}). Finally, the user has to define a parameter of type \texttt{int} that counts the number of Leja iterations per time step. This may be used as a proxy for the overall computational runtime - the larger the number of iterations, the higher the runtime. As \texttt{real\_Leja\_phi\_nl} (Listing \ref{alg:leja_phi_nl_ext}) is to be used to evaluate nonhomogeous linear differential equations, it requires one additional parameter: the $\varphi$-function that is to be applied to \texttt{u\_input}. This depends on the specific exponential method used. A demonstration of calling the functions \texttt{real\_Leja\_exp} and \texttt{real\_Leja\_phi\_nl} can be found in \lexint $\rightarrow$ \texttt{Test} $\rightarrow$ \texttt{test\_2D.cu} and \lexint $\rightarrow$ \texttt{Test} $\rightarrow$ \texttt{test\_2D.cpp}.

\subsubsection*{Invoking the exponential integrators}

\begin{minipage}{0.95\linewidth}
\begin{lstlisting}[language = C++, caption = Embedded exponential integrator, label = alg:exp_int_gpu]
//? Solvers without an error estimate
template <typename rhs>
void Leja_GPU::exp_int(rhs& RHS,               //? RHS function
                       double* u_input,        //? Input state variable(s)
                       double* u_output,       //? Output
                       double c,               //? Shifting factor
                       double Gamma,           //? Scaling factor
                       double rtol,            //? Relative tolerance (normalised)
                       double atol,            //? Absolute tolerance
                       double dt,              //? Step size
                       int& iters,             //? Number of Leja iterations
                       bool GPU                //? false --> CPU; true --> GPU
                       )
\end{lstlisting}

\begin{lstlisting}[language = C++, caption = Embedded exponential integrator, label = alg:embed_exp_int_gpu]
//? Embedded integrators
template <typename rhs>
void Leja_GPU::embed_exp_int(rhs& RHS,               //? RHS function
                             double* u_input,        //? Input state variable(s)
                             double* u_output_low,   //? Output (lower order)
                             double* u_output_high,  //? Output (higher order)
                             double& error,          //? Error estimate
                             double c,               //? Shifting factor
                             double Gamma,           //? Scaling factor
                             double rtol,            //? Relative tolerance(normalised)
                             double atol,            //? Absolute tolerance
                             double dt,              //? Step size
                             int& iters,             //? Number of Leja iterations
                             bool GPU                //? false --> CPU; true --> GPU
                             )
\end{lstlisting}
\end{minipage}

The class \texttt{Leja} has two member functions to invoke the desired exponential integrator, namely \texttt{exp\_int} and \texttt{embed\_exp\_int}. \texttt{exp\_int} takes the user-defined \texttt{RHS} function, the input vector (\texttt{u\_input}), and yields the output vector (\texttt{u\_output}) for integrators without an embedded error estimate, whilst \texttt{embed\_exp\_int} yields two output vectors - \texttt{u\_output\_low} and \texttt{u\_output\_high}, corresponding to the lower-order and the higher-order solutions, respectively, for exponential integrators with an embedded error estimate (see Listings \ref{alg:exp_int_gpu} and \ref{alg:embed_exp_int_gpu} and the descriptions therein). The difference between these two vectors yields an error estimate, the l2 norm normalised to $\sqrt{N}$ of which is returned as \texttt{error}. This parameter may be used to control the step sizes in case of a variable step size implementation. 

Other arguments of these functions are similar to \texttt{real\_Leja\_exp} and \texttt{real\_Leja\_phi\_nl} - these include the Leja interpolation parameters (\texttt{c} and \texttt{Gamma}), user-defined tolerances (\texttt{rtol} and \texttt{atol}), step size (\texttt{dt}) and the specification as to whether GPU support is enabled. These parameters are to be specified by the user whilst invoking an exponential integrator. An outline of how an embedded exponential integrator should be called is shown in Listing \ref{alg:lexint}, respectively. Further details and related technical aspects can be found in the test examples - \lexint $\rightarrow$ \texttt{Test} $\rightarrow$ \texttt{test\_2D.cu} and \lexint $\rightarrow$ \texttt{Test} $\rightarrow$ \texttt{test\_2D.cpp}.

\begin{minipage}{0.95\linewidth}
\begin{lstlisting}[language = C++, caption = Calling an exponential integrator using \lexint, label = alg:lexint]
#include "./LeXInt/CUDA/Leja.hpp"           // Check for proper directory

//? Define the desired RHS function
rhs RHS(*args)                              //* RHS of desired problem

//? Create an object to call LeXInt functions
Leja<rhs> leja_gpu{N, desired_integrator = "EPIRK4s3A"};   

//? Allocate memory on GPU for input and output vectors
size_t N_size = N * sizeof(double);
double *device_u; cudaMalloc(&device_u, N_size);              //* Input vector
double *device_u_sol; cudaMalloc(&device_u_sol, N_size);      //* Higher order solution
double *device_u_low; cudaMalloc(&device_u_low, N_size);      //* Lower order solution
double error;                                                 //* Error estimate
int iters;                                                    //* # of Leja iterations

//? Compute spectrum (shifting and scaling parameters)
double eigenvalue;
leja_gpu.Power_iterations(RHS, device_u, eigenvalue, GPU = true);
eigenvalue = -1.05*eigenvalue;              // Real eigenvalue has to be negative
double c = eigenvalue/2.0; double Gamma = -eigenvalue/4.0;

//? Time loop
while (time < t_final)
{
    //? Error embedded exponential integrator
    leja_gpu.embed_exp_int(RHS, device_u, device_u_low, device_u_sol, 
                           error, c, Gamma, tol, dt, iters, GPU = true);


    //? ** Recompute spectrum (c and Gamma) every N time steps **
    leja_gpu.Power_iterations(RHS, device_u, eigenvalue, GPU = true);
    eigenvalue = -1.05*eigenvalue;              // Real eigenvalue has to be negative
    double c = eigenvalue/2.0; double Gamma = -eigenvalue/4.0;

    //* Update variables
    time = time + dt;
    LeXInt::copy(device_u_sol, device_u, N, GPU = true);
}
cudaDeviceSynchronize();
\end{lstlisting}
\end{minipage}

%%% --------------------------------------------------------- %%%

\subsection{Code structure (``internal")}

In this subsection, we provide a low-level description of some of the technical aspects of \lexint. The struct \texttt{Leja}, in \texttt{Leja.hpp}, stores a range of helper functions and handles the memory management for the auxiliary variables required for each of the exponential integrators and the functions used to compute the interpolation on Leja points and the action of the Jacobian on the relevant vector.

\subsubsection*{Leja interpolation}
Listings \ref{alg:leja_exp_int} and \ref{alg:leja_phi_nl_int} illustrate the ``internal'' structure of the functions \texttt{real\_Leja\_exp} and \texttt{real\_Leja\_phi\_nl}, respectively. In \texttt{real\_Leja\_exp}, \texttt{u} corresponds to the input state variable or the set of state variables and \texttt{polynomial} is the resulting output interpolated on Leja points. \texttt{auxiliary\_Leja} (see Listing \ref{alg:internal_arrays}) stores four additional `N'-dimensional arrays that are needed for computing the matrix exponential (or $\varphi_l$ function) applied to a vector and approximating the spectrum of the underlying matrix, i.e, to compute the most dominant eigenvalue. \texttt{real\_Leja\_exp} and \texttt{real\_Leja\_phi\_nl} need only one of these four additional `N'-dimensional arrays, whereas all four are needed to approximate the spectrum. To avoid allocating and deallocating memory on the GPU at every time step, this memory is allocated only once - when an object of this class is created. Two additional parameters of the ``internal'' function are \texttt{Leja\_X} (the array of Leja points) and \texttt{N} (the number of grid points). The parameters that are exclusive to the ``internal" function (\texttt{auxiliary\_Leja}, \texttt{Leja\_X}, and \texttt{N}) are automatically taken care of by \lexint.

The only difference in \texttt{real\_Leja\_phi\_nl} (Listing \ref{alg:leja_phi_nl_int}) is that the input array, that is to be interpolated on Leja points, is named \texttt{interp\_vector}. This, often times, is different from the state variable \texttt{u}, for nonhomogenous differential equations.

\begin{minipage}{0.95\linewidth}
\begin{lstlisting}[language = C++, caption = Leja interpolation of the matrix exponential applied to vector `u' (``internal" function), label = alg:leja_exp_int]
//? Real Leja Exp ( Homogenous linear equations - matrix exponential )
template <typename rhs>
void real_Leja_exp(rhs& RHS,                   //? RHS function
                   double* u,                  //? Input state variable(s)
                   double* polynomial,         //? Output matrix exponential times 'u'
                   double* auxiliary_Leja,     //? Internal auxiliary variables (Leja)
                   size_t N,                   //? Number of grid points
                   vector<double>& Leja_X,     //? Array of Leja points
                   double c,                   //? Shifting factor
                   double Gamma,               //? Scaling factor
                   double rtol,                //? Relative tolerance (normalised)
                   double atol,                //? Absolute tolerance
                   double dt,                  //? Step size
                   int& iters,                 //? # of iterations needed to converge 
                   bool GPU,                   //? false (0) --> CPU; true (1) --> GPU
                   GPU_handle& cublas_handle   //? CuBLAS handle
                   )
\end{lstlisting}

\begin{lstlisting}[language = C++, caption = Leja interpolation of the $\varphi$ function applied to \texttt{interp\_vector} (``internal" function)., label = alg:leja_phi_nl_int]
//? Real Leja Phi NL (Nonhomogenous linear equations)
template <typename rhs>
void real_Leja_phi_nl(
        rhs& RHS,                           //? RHS function
        double* interp_vector,              //? Input state variable(s)
        double* polynomial,                 //? Output applied to phi function
        double* auxiliary_Leja,             //? Internal auxiliary variables (Leja)
        size_t N,                           //? Number of grid points
        double (* phi_function) (double),   //? Phi function
        vector<double>& Leja_X,             //? Array of Leja points
        double c,                           //? Shifting factor
        double Gamma,                       //? Scaling factor
        double rtol,                        //? Relative tolerance (normalised)
        double atol,                        //? Absolute tolerance
        double dt,                          //? Step size
        int& iters,                         //? # of iterations needed to converge
        bool GPU,                           //? false (0) --> CPU; true (1) --> GPU
        GPU_handle& cublas_handle           //? CuBLAS handle
        )
\end{lstlisting}
\end{minipage}

\texttt{real\_Leja\_phi} (Listing \ref{alg:leja_phi}) is used to interpolate the $\varphi_l(z)$ appearing in the exponential integrators. This function is similar to \texttt{real\_Leja\_phi\_nl} to a great extent. \texttt{u} is the input state variable needed to compute the action of the Jacobian on the vector to be applied to \texttt{phi\_function} interpolated on Leja points, the output of which is stored in \texttt{polynomial}. \texttt{real\_Leja\_phi} needs a variable to indicate the number of vertical \cite{Tokman16} interpolations (\texttt{integrator\_coeffs}) that is specific to the integrator under consideration. This is handled implicitly in \lexint. We note that whilst \texttt{real\_Leja\_exp} and \texttt{real\_Leja\_phi\_nl} need only one `N'-dimensional array (stored in \texttt{auxiliary\_Leja}), \texttt{real\_Leja\_phi} needs four of the same. These auxiliary arrays are used to compute the action of the Jacobian on the relevant vector, nonlinear remainders in exponential integrators, and to compute the most dominant eigenvalue using power iterations.

\subsubsection*{Invoking the exponential integrators}
\texttt{auxiliary\_expint} stores memory for a certain number of internal arrays, which we call `\texttt{num\_vectors}' (Listing \ref{alg:internal_arrays}). This number depends on the exponential integrator under consideration. Each of these vectors consists of `\texttt{N}' \texttt{doubles}, where `\texttt{N}' is the total number of grid points. \texttt{auxiliary\_expint} is automatically passed as function arguments to the exponential integrators. The memory allocated for \texttt{auxiliary\_expint} and \texttt{auxiliary\_Leja} is freed at the end of the computations.

\begin{minipage}{0.95\linewidth}
\begin{lstlisting}[language = C++, caption = Allocation of memory for the intenal arrays for Leja interoplation and exponential integrators., label = alg:internal_arrays]
//! Allocate memory - these are device vectors if GPU support is activated
double* auxiliary_Leja;    //? Internal vectors for Leja interpolation and power iters
double* auxiliary_expint;  //? Internal vectors for an exponential integrator

//! Constructor
Leja(int _N, string _integrator_name) :  N(_N), integrator_name(_integrator_name)
{
    if (integrator_name == "Rosenbrock_Euler")
    {
        num_vectors = 1;
    }
    else if (integrator_name == "EXPRB32")
    {
        num_vectors = 1;
    }
    else if (integrator_name == "EXPRB42")
    {
        num_vectors = 2;
    }
    .
    .
    .
}

#ifdef __CUDACC__
    //? Allocate memory on device
    cudaMalloc(&auxiliary_Leja, 4 * N * sizeof(double));
    cudaMalloc(&auxiliary_expint, num_vectors * N * sizeof(double));
#else
    auxiliary_Leja = (double*)malloc(4 * N * sizeof(double));
    auxiliary_expint = (double*)malloc(num_vectors * N * sizeof(double));
#endif

//! Destructor
~Leja()
{
    //? Deallocate memory
    #ifdef __CUDACC__
        cudaFree(auxiliary_Leja);
        cudaFree(auxiliary_expint);        
    #else
        free(auxiliary_Leja);
        free(auxiliary_expint);
    #endif
}
\end{lstlisting}
\end{minipage}

Let us consider the examples of two exponential integrators, Rosenbrock--Euler (Listing \ref{alg:Ros_Eu}) and EXPRB32 (Listing \ref{alg:exprb32}), the equations of which are given by
\[u^{n+1} = u^n + \varphi_1\left(\mathcal{J}(u^n) \Delta t\right) f(u^n) \Delta t\]
and 
\[a^n = u^n + \varphi_1\left(\mathcal{J}(u^n) \Delta t\right) f(u^n) \Delta t, \qquad
u^{n+1} = a^n + 2\varphi_3\left(\mathcal{J}(u^n) \Delta t\right) (\mathcal{F}(a^n) - \mathcal{F}(u^n)) \Delta t,\]
respectively. Here, $\mathcal{F}(k^n) = f(k^n) - \mathcal{J}(u^n)k^n$ and $\mathcal{J}(u^n)$ is the Jacobian matrix evaluated at $u^n$. The action of the Jacobian on the relevant vector is computed numerically using finite differences. We need only one internal vector for the Rosenbrock--Euler integrator, which we call \texttt{f\_u}. The RHS function evaluated at $u^n$ multiplied by the time step size, i.e. $f(u^n) \Delta t$, is stored in \texttt{f\_u}. The action of the $\varphi_l$ function on the relevant vector, here $\varphi_1\left(\mathcal{J}(u^n) \Delta t\right) f(u^n) \Delta t$, which is the output of the function \texttt{real\_Leja\_phi}, is stored in the output vector \texttt{u\_exprb2}. The final solution is obtained by adding \texttt{u} and \texttt{u\_exprb2}, which is subsequently overwritten to \texttt{u\_exprb2}.

EXPRB32 is an embedded exponential integrator, and as such, we now have two output vectors. If these output vectors are used sparingly in the internal stages, the number of additional vectors needed can be significantly reduced. In EXPRB32, we need only one additional vector. We name the internal vectors in a way that improves the readability of the code. The names of these vectors are similar to the ones in the \texttt{Python} version of \lexint. Similar to the Rosenbrock--Euler integrator, in EXPRB32 (Listing \ref{alg:exprb32}), we first compute $f(u^n) \Delta t$ and store the output vector in the auxiliary vector that we call \texttt{f\_u}.  The action of the $\varphi_l$ function on the relevant vector, here $\varphi_1\left(\mathcal{J}(u^n) \Delta t\right) f(u^n) \Delta t$, which is the output of the function \texttt{real\_Leja\_phi} is stored in \texttt{u\_exprb2}, after which, we add \texttt{u} and \texttt{u\_exprb2} to get the second-order solution. Now, we compute the nonlinear remainders $\mathcal{F}(a^n) \Delta t$ and $\mathcal{F}(u^n) \Delta t$ and store them in \texttt{auxiliary\_expint} and \texttt{u\_exprb3}, respectively. We call these vectors \texttt{NL\_u} and \texttt{NL\_a}, respectively. The difference of these two vectors, called \texttt{R\_a}, is overwritten to the memory of \texttt{u\_exprb3}. We then approximate $\varphi_3\left(\mathcal{J}(u^n) \Delta t\right) (\mathcal{F}(a^n) - \mathcal{F}(u^n)) \Delta t$ and the resulting polynomial is stored in \texttt{auxiliary\_expint}, which we now call \texttt{u\_nl\_3}. Finally, the third-order solution, \texttt{u\_exprb3}, is obtained by summing up \texttt{u\_exprb2} and $2 \cdot$ \texttt{u\_nl\_3}. The error estimate is obtained by computing the norm of $2 \cdot$ \texttt{u\_nl\_3}. In a similar fashion, we have accounted for a highly optimised number of internal vectors for all integrators in \texttt{Leja.hpp}. It is worth noting that higher-order integrators may require a large number of internal vectors, which translates to hefty memory requirements, thereby potentially rendering them infeasible for practical applications.

\begin{minipage}{0.95\linewidth}
\begin{lstlisting}[language = C++, caption = Leja interpolation of the $\varphi$ function applied to vector \texttt{interp\_vector} for the exponential integrators (``internal" function)., label = alg:leja_phi]
//? Real Leja Phi (Exponential integrators for nonlinear equations)
template <typename rhs>
void real_Leja_phi(
        rhs& RHS,                           //? RHS function
        double* u,                          //? Input state variable(s)
        double* interp_vector,              //? Input to be multiplied to phi function
        double* polynomial,                 //? Output multiplied to phi function
        double* auxiliary_Leja,             //? Internal auxiliary variables (Leja)
        size_t N,                           //? Number of grid points
        vector<double> integrator_coeffs,   //? Coefficients of the integrator
        double (* phi_function) (double),   //? Phi function
        vector<double>& Leja_X,             //? Array of Leja points
        double c,                           //? Shifting factor
        double Gamma,                       //? Scaling factor
        double rtol,                        //? Relative tolerance (normalised)
        double atol,                        //? Absolute tolerance
        double dt,                          //? Step size
        int& iters,                         //? # of iterations needed to converge
        bool GPU,                           //? false (0) --> CPU; true (1) --> GPU
        GPU_handle& cublas_handle           //? CuBLAS handle
        )
\end{lstlisting}

\begin{lstlisting}[language = C++, caption = Rosenbrock--Euler, label = alg:Ros_Eu]
template <typename rhs>
void Ros_Eu(rhs& RHS,                   //? RHS function
            double* u,                  //? Input state variable(s)
            double* u_exprb2,           //? Output state variable(s)
            double* auxiliary_expint,   //? Internal auxiliary variables
            double* auxiliary_Leja,     //? Internal auxiliary variables (Leja)
            size_t N,                   //? Number of grid points
            vector<double>& Leja_X,     //? Array of Leja points
            double c,                   //? Shifting factor
            double Gamma,               //? Scaling factor
            double rtol,                //? Relative tolerance (normalised)
            double atol,                //? Absolute tolerance
            double dt,                  //? Step size
            int& iters,                 //? # of iterations needed to converge
            bool GPU,                   //? false (0) --> CPU; true (1) --> GPU
            GPU_handle& cublas_handle   //? CuBLAS handle
            )
{
    //? Assign names and variables
    double* f_u = &auxiliary_expint[0];

    //? RHS evaluated at 'u' multiplied by 'dt'; f_u = RHS(u)*dt
    RHS(u, f_u);
    axpby(dt, f_u, f_u, N, GPU);

    //? Interpolation of RHS(u) at 1; phi_1(J(u) dt) f(u) dt
    real_Leja_phi(RHS, u, f_u, u_exprb2, auxiliary_Leja, N, {1.0}, 
                    phi_1, Leja_X, c, Gamma, tol, dt, iters, GPU, cublas_handle);

    //? 2nd order solution; u_2 = u + phi_1(J(u) dt) f(u) dt
    axpby(1.0, u, 1.0, u_exprb2, u_exprb2, N, GPU);
}
\end{lstlisting}
\end{minipage}

\begin{minipage}{0.95\linewidth}
\begin{lstlisting}[language = C++, caption = EXPRB32, label = alg:exprb32]
template <typename rhs>
void EXPRB32(rhs& RHS,                   //? RHS function
             double* u,                  //? Input state variable(s)
             double* u_exprb2,           //? Output state variable(s) (lower order)
             double* u_exprb3,           //? Output state variable(s) (higher order)
             double& error,              //? Embedded error estimate
             double* auxiliary_expint,   //? Auxiliary variables (EXPRB32)
             double* auxiliary_Leja,     //? Auxiliary variables (Leja and NL arrays)
             size_t N,                   //? Number of grid points
             vector<double>& Leja_X,     //? Array of Leja points
             double c,                   //? Shifting factor
             double Gamma,               //? Scaling factor
             double rtol,                //? Relative tolerance (normalised)
             double atol,                //? Absolute tolerance
             double dt,                  //? Step size
             int& iters,                 //? # of iterations needed to converge
             bool GPU,                   //? false (0) --> CPU; true (1) --> GPU
             GPU_handle& cublas_handle   //? CuBLAS handle
             )
{
    //? Counters for Leja iterations
    int iters_1 = 0, iters_2 = 0;

    //? Assign names and variables
    double* f_u = &auxiliary_expint[0]; double* u_flux = &u_exprb2[0]; 
    double* NL_u = &auxiliary_expint[0]; double* NL_a = &u_exprb3[0]; 
    double* R_a = &u_exprb3[0]; 
    double* u_nl_3 = &auxiliary_expint[0]; double* error_vector = &auxiliary_expint[0];

    //? RHS evaluated at 'u' multiplied by 'dt'; f_u = RHS(u)*dt
    RHS(u, f_u);
    axpby(dt, f_u, f_u, N, GPU);

    //? Interpolation of RHS(u) at 1; u_flux = phi_1(J(u) dt) f(u) dt
    real_Leja_phi(RHS, u, f_u, u_flux, auxiliary_Leja, N, {1.0}, 
                    phi_1, Leja_X, c, Gamma, tol, dt, iters_1, GPU, cublas_handle);

    //! Internal stage 1; 2nd order solution; u_2 = u + phi_1(J(u) dt) f(u) dt
    axpby(1.0, u, 1.0, u_flux, u_exprb2, N, GPU);

    //? R_a = (NL_a - NL_u) * dt
    Nonlinear_remainder(RHS, u, u,        NL_u, auxiliary_Leja, N, GPU, cublas_handle);
    Nonlinear_remainder(RHS, u, u_exprb2, NL_a, auxiliary_Leja, N, GPU, cublas_handle);
    axpby(dt, NL_a, -dt, NL_u, R_a, N, GPU);

    //? u_nl_3 = phi_3(J(u) dt) R(a) dt
    real_Leja_phi(RHS, u, R_a, u_nl_3, auxiliary_Leja, N, {1.0}, 
                    phi_3, Leja_X, c, Gamma, tol, dt, iters_2, GPU, cublas_handle);
                    
    //! 3rd order solution; u_3 = u_2 + 2 phi_3(J(u) dt) R(a) dt
    axpby(1.0, u_exprb2, 2.0, u_nl_3, u_exprb3, N, GPU);

    //? Error estimate
    axpby(2.0, u_nl_3, error_vector, N, GPU);
    error = l2norm(error_vector, N, GPU, cublas_handle);

    //? Total number of Leja iterations
    iters = iters_1 + iters_2;
}
\end{lstlisting}
\end{minipage}

%%%%%%%%%%%%%%%%%%%%%%%%%%%%%%%%%%%%%%%%%%%%%%%%%%%%%%%%%%%%%%%%%%%%%%%%%%%%%%%%%%

\section{Illustrative examples}
\label{sec:eg}

We compare and contrast the performance of the \texttt{C++} and \texttt{CUDA} codes for a range of two-dimensional test problems. We consider periodic boundary conditions on $x \times y \in [-1, 1] \times [-1, 1]$ in all cases and $\nu$ corresponds to the advection velocity, wherever applicable. We note that the operators, $\nabla$ and $\nabla^2$, are discretised using the third-order upwind and second-order centred difference schemes, respectively.  

The problems we are interested in this work typically result from the discretisation of ODEs/PDEs. This results in either a stencil code or structured and sparse matrix-vector products that are memory-bound problems. Examples include advection--diffusion--reaction equations, Navier--Stokes equation, the set of equations for magnetohydrodynamics, the Maxwell-Vlasov equation, etc. For the memory-bound test problems considered in this work, we use the memory bandwidth as a metric to study the performance of \lexint. We compute the average bandwidth achieved (on the A100 GPU) as
\[ \text{Bandwidth} = \frac{N_\mathrm{grid} \times \texttt{sizeof(double)} \times N_\mathrm{rw} \times 10^{-9} }{\mathrm{Simulation \; time}}  \, \mathrm{GB/s}, \]
where $N_\mathrm{rw}$ is the total number of vector reads and writes (of size $N_\mathrm{grid}$) throughout the simulations. Let us note that runtime is an absolute measure of performance (and this is what practitioners care about), the bandwidth achieved is an ``relative'' measure - it indicates how close one is to achieving the maximum possible throughput on a certain device (for memory-bound problems). Other discretisation methods such as finite element or discontinuous Galerkin may tend to be algorithmically intensive, and one might be able to leverage GPUs to speedup these computations. This is beyond the scope of the current work.

The GPU simulations were carried out on a single A100 GPU on the Booster partition, whereas the CPU simulations were carried out on the Data-Centric General Partition (DCGP) partition of the Leonardo cluster \footnote{\url{https://www.hpc.cineca.it/systems/hardware/leonardo/}} hosted at CINECA \footnote{\url{https://leonardo-supercomputer.cineca.eu}}. The theoretical maximum achievable bandwidth on an A100 is 1555 GB/s, however, this is rarely practically achievable. Even the most elementary \texttt{CUDA} functions such as \texttt{cudaMemcpy} (device to device) or vector operations such as $z = ax + by$ fail to run with the theoretical maximum bandwidth. Complex operations are obviously expected to exhibit lower bandwidths. For the two aforementioned simple test problems on an A100 Leonardo Booster, \texttt{cudaMemcpy} (device to device) and $z = ax + by$, we achieved bandwidths of 1359.14 GB/s and 1438.48 GB/s, respectively, with proper care taken to saturate the memory on the entire node. All bandwidths reported in this article are normalised to the maximum bandwidth achieved for the simple test problem of vector operations, i.e., 1438.48 GB/s. 

\textbf{Problem I: Diffusion--Advection} - The homogeneous linear diffusion--advection equation reads
\begin{equation}
    \frac{\partial u}{\partial t} = \frac{\partial^2 u}{\partial x^2} + \frac{\partial^2 u}{\partial y^2} + \nu \, \left(\frac{\partial u}{\partial x} + \frac{\partial u}{\partial y}\right).
    \label{eq:da}
\end{equation}
We choose $\nu = 10$ and the initial condition to be $u(x, y, t = 0) = 1 + \exp(-((x + 0.5)^2 + (y + 0.5)^2)/0.01)$. Eq. \eqref{eq:da} can be written as \[\frac{\partial u}{\partial t} = \mathcal{A}_1 u,\] where $\mathcal{A}_1 = \nabla^2 (\cdot) + \nu \, \nabla (\cdot)$. The solution to this equation is given by $u^{n+1} = \exp(\mathcal{A}_1 \Delta t) u^n$. Integrating this equation in time using \lexint requires invoking the \texttt{LeXInt::real\_Leja\_exp} function that returns the action of the matrix exponential on the desired vector, i.e. $\exp(\mathcal{A}_1 \Delta t) u^n$.

We consider two different grid resolutions ($2^{13} \times 2^{13}$ and $2^{14} \times 2^{14}$) with different simulation times. The average normalised bandwidth achieved on the A100 GPU, for problem I, is reported in Table \ref{tab:das} and the speedup achieved on the GPU over the CPU implementation using 1 and 112 core(s) on a dual-socket Intel Xeon Platinum 8480p CPU system on Leonardo DCGP, for selected parameters, is illustrated in Table \ref{tab:speedups}.

\begin{figure}
    \centering
    \includegraphics[width = 0.495\columnwidth]{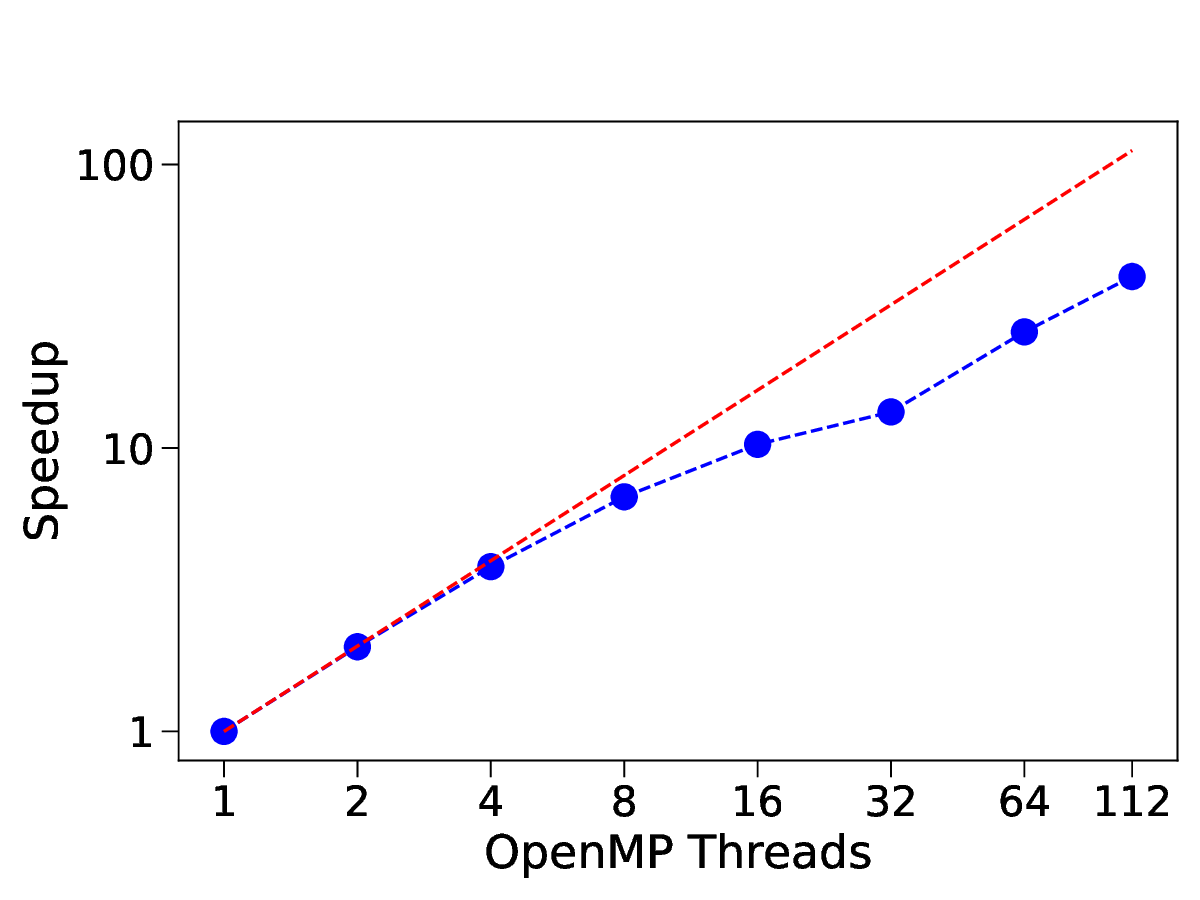}
    \includegraphics[width = 0.495\columnwidth]{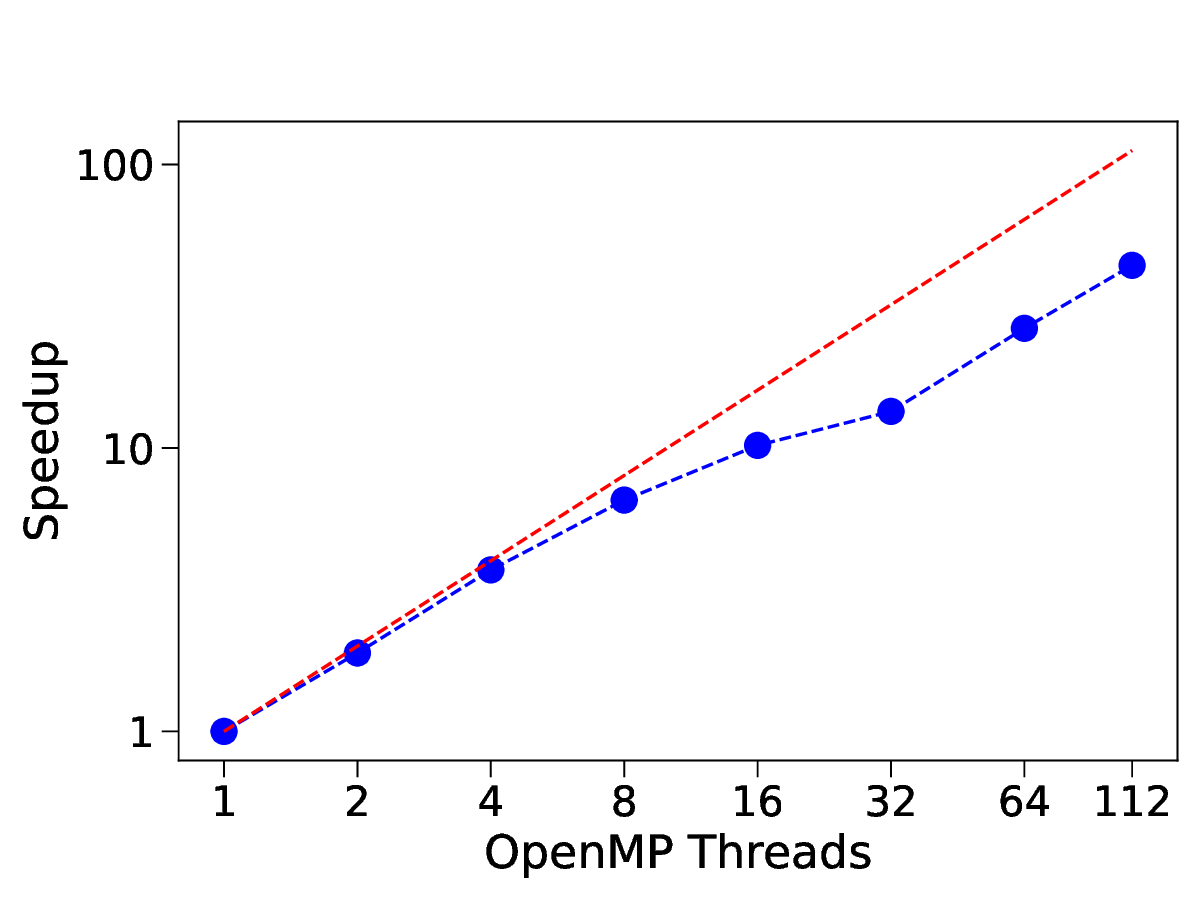} \\
    \includegraphics[width = 0.495\columnwidth]{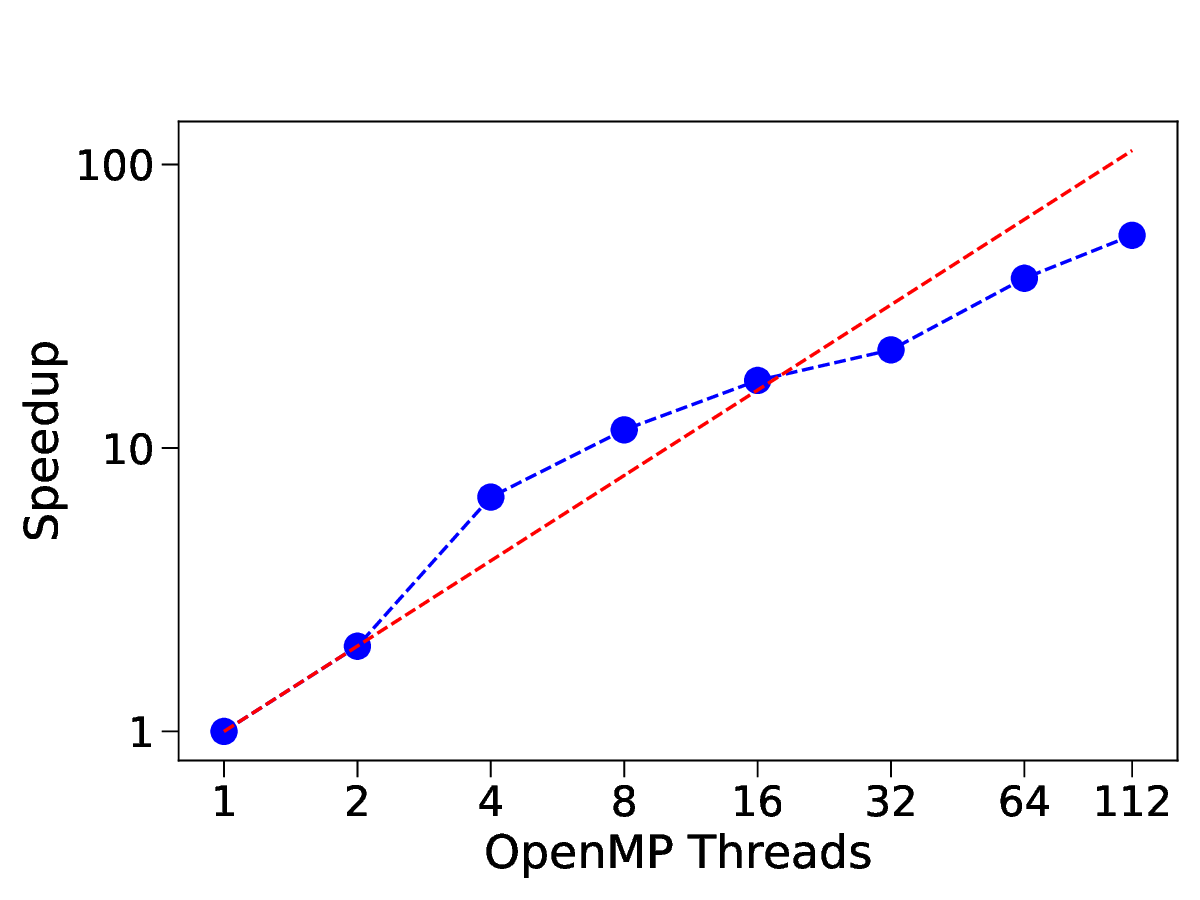}
    \includegraphics[width = 0.495\columnwidth]{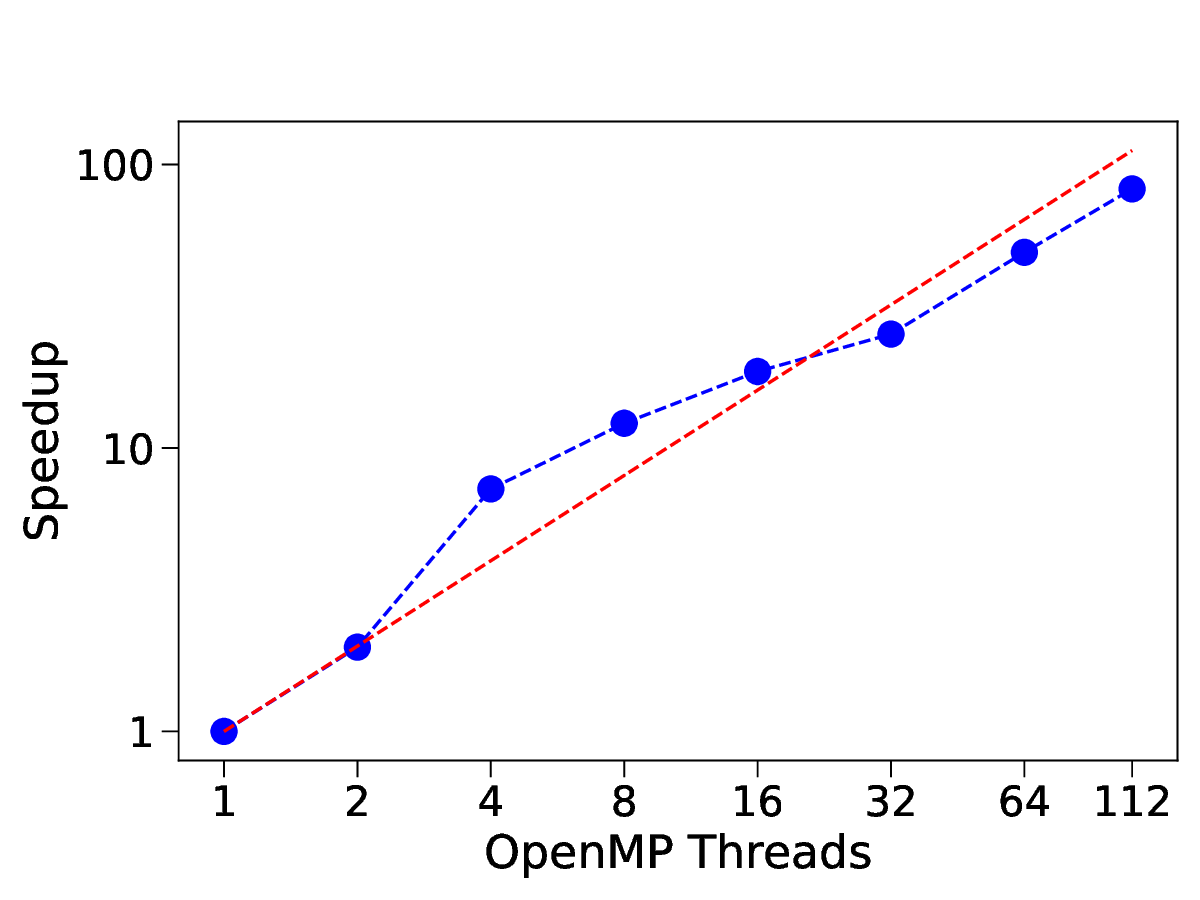} \\
    \includegraphics[width = 0.495\columnwidth]{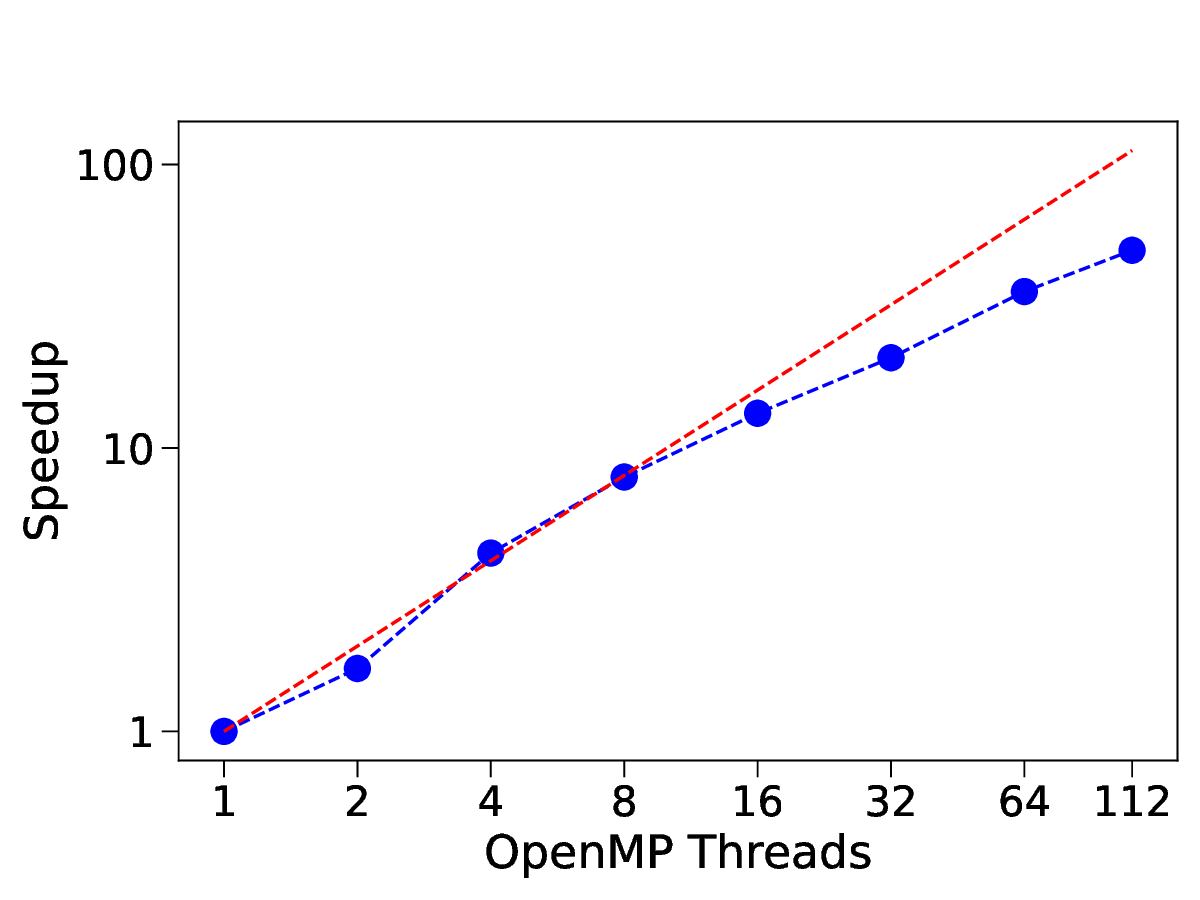}
    \includegraphics[width = 0.495\columnwidth]{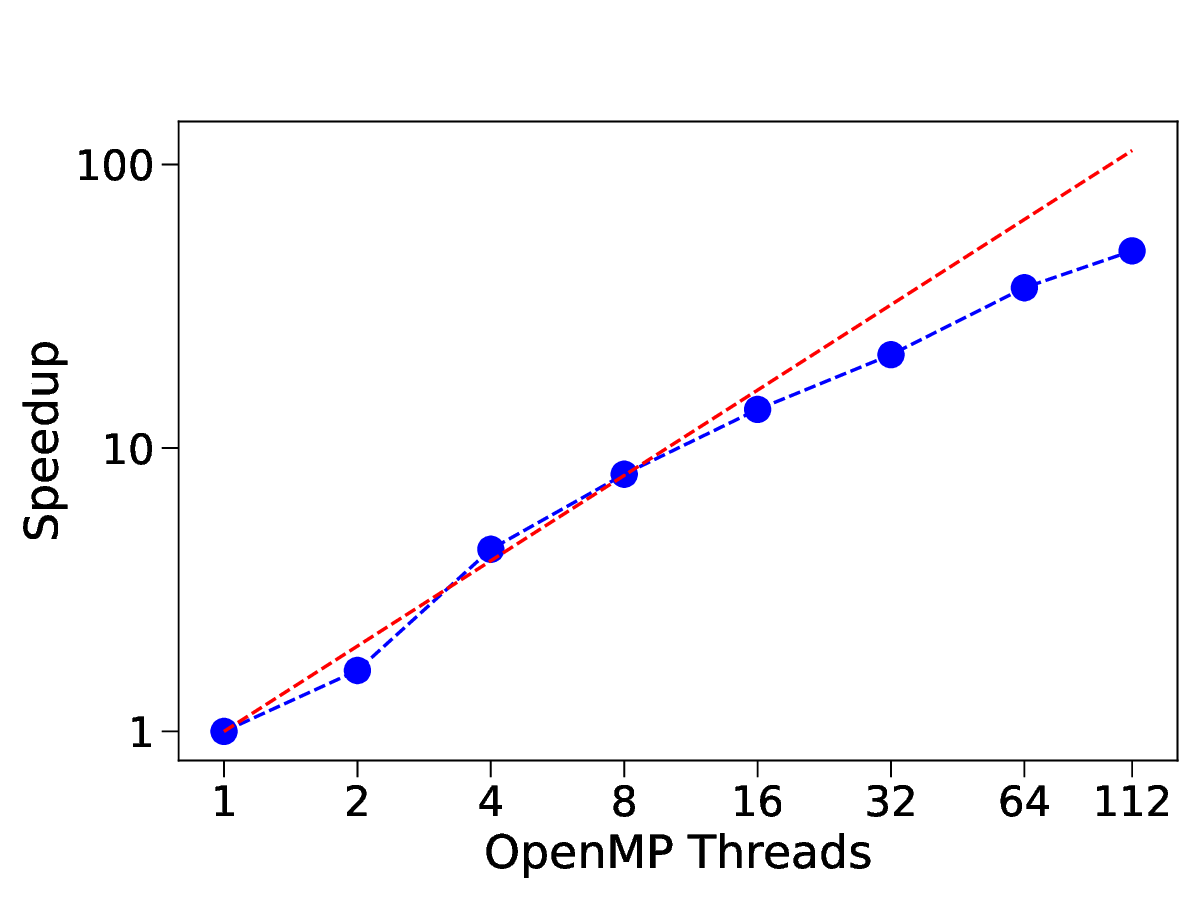}
    \caption{The blue circles illustrate strong scaling for the linear diffusion--advection problem (top panel), linear diffusion--advection problem in the presence of sources (middle panel), and the Burgers' equation integrated using EXPRB32 (bottom panel). The configuration for simulations for the plots on the left half is $N_\mathrm{grid} = 2^{13} \times 2^{13}$, $T_f = 2\cdot10^{-7}$, and $\Delta t = 1 \cdot \Delta t_\mathrm{CFL}$, whilst for the right half, it is $N_\mathrm{grid} = 2^{14} \times 2^{14}$, $T_f = 2\cdot10^{-6}$, and $\Delta t = 100 \cdot \Delta t_\mathrm{CFL}$. All simulations were performed with a user-specified tolerance of $10^{-12}$ on Leonardo DCGP. The red line represents ideal strong scaling.}
    \label{fig:strong_scaling}
\end{figure}

We note that a CPU-GPU performance comparison is worthy only if the code has been well optimised on CPUs. An ill-optimised code would show severe degradation in a strong scaling test. Here, we aim to show that the GPU version has significantly better performance than a well-optimised CPU code. As such, we present strong scaling tests for the selected parameters in Fig. \ref{fig:strong_scaling} for the CPU version of the code. The simulations for the strong scaling were performed on a single node of Leonardo DCGP that has a maximum of 112 cores per node.

\textbf{Problem II: Diffusion--Advection (with source)} - The nonhomogeneous linear diffusion--advection equation in the presence of time-independent sources is given by 
\begin{equation}
    \frac{\partial u}{\partial t} = \frac{\partial^2 u}{\partial x^2} + \frac{\partial^2 u}{\partial y^2} + \nu \, \left(\frac{\partial u}{\partial x} + \frac{\partial u}{\partial y}\right) + S(x, y).
    \label{eq:ds}
\end{equation}
The initial condition is chosen to be exactly the same as in the previous example and $\nu = 10$. We choose the source as $S(x, y) = \exp(-((x + 0.4)^2 + (y - 0.6)^2)/0.05) + \exp(-((x - 0.25)^2 + (y + 0.1)^2)/0.04)$). Rewriting Eq. \eqref{eq:ds} as \[\frac{\partial u}{\partial t} = \mathcal{A}_1 u + S(x, y),\] where $\mathcal{A}_1 = \nabla^2 (\cdot) + \nu \, \nabla (\cdot)$, the solution to this equation is $u^{n+1} = u^n + \varphi_1(\mathcal{A}_1 \Delta t) (\mathcal{A}_1 u^n + S(x, y))$. Integrating this equation in time using \lexint requires invoking the \texttt{LeXInt::real\_Leja\_phi\_nl} function that returns $\varphi_l(\mathcal{A}_1 \Delta t) (\mathcal{A}_1 u^n + S(x, y))$. Similar to the previous problem, two different grid resolutions are considered (with different simulation times): $2^{13} \times 2^{13}$ and $2^{14} \times 2^{14}$, and resulting average normalised bandwidth achieved on the A100 GPU are summarised in Table \ref{tab:das}. 

\textbf{Problem III: Burgers' Equation} - The nonlinear Burgers' equation, in the presence of viscosity, reads 
\begin{equation}
   \frac{\partial u}{\partial t} = \frac{\partial^2 u}{\partial x^2} + \frac{\partial^2 u}{\partial y^2} + \frac{\nu}{2} \, \left(\frac{\partial u^2}{\partial x} + \frac{\partial u^2}{\partial y}\right).
    \label{eq:burgers}
\end{equation}
Here, $\nu = 10$ and the initial condition is chosen to be \[u(x, y, t = 0) = 2 + 10^{-2} \left[\sin(2 \pi x) + \sin(2 \pi y) + \sin(8 \pi x + 0.3) + \sin(8 \pi y + 0.3)\right]. \]
In Table \ref{tab:bur}, we show the bandwidth achieved on the A100 for second-order Rosenbrock--Euler, third-order EXPRB32, fourth-order EPIRK4s3A, and fifth-order EPIRK5P1 solvers for different grid resolutions and step sizes. Owing to the memory requirements of high-order integrators, a resolution of $2^{14} \times 2^{14}$ is not feasible for EPIRK5P1. The use of exponential integrators involves invoking the \texttt{LeXInt::real\_Leja\_phi} function that computes the $\varphi_l$ function on the relevant vector. \texttt{LeXInt::real\_Leja\_phi} is designed to facilitate computation of the $\varphi_l$ functions in ``vertical'' \cite{Tokman16} wherever possible.

\begin{table}[t]
    \centering
    \begin{tabular}{|c|c|c|c|c|} 
        \hline
        $N_\mathrm{grid}$ & $T_f$ & $\Delta t$ & Bandwidth & Bandwidth \\
        &  &  & (Problem I) & (Problem II) \\
        \hline 
                               & $2\cdot10^{-7}$ & $1   \cdot \Delta t_\mathrm{CFL}$  & 0.87 & 0.89 \\ 
        $2^{13} \times 2^{13}$ & $2\cdot10^{-6}$ & $10  \cdot \Delta t_\mathrm{CFL}$  & 0.89 & 0.90 \\
                               & $2\cdot10^{-5}$ & $100 \cdot \Delta t_\mathrm{CFL}$  & 0.91 & 0.90 \\ 
        \hline
                               & $2\cdot10^{-8}$ & $1   \cdot \Delta t_\mathrm{CFL}$  & 0.89 & 0.90 \\ 
        $2^{14} \times 2^{14}$ & $2\cdot10^{-7}$ & $10  \cdot \Delta t_\mathrm{CFL}$  & 0.91 & 0.92 \\ 
                               & $2\cdot10^{-6}$ & $100 \cdot \Delta t_\mathrm{CFL}$  & 0.92 & 0.92 \\ 
        \hline
    \end{tabular}
    \caption{Bandwidth achieved on A100, normalised to the maximum practically-achieved bandwidth, for problems I and II with different grid resolutions and time step sizes. All simulations were performed with a user-specified tolerance of $10^{-12}$.}
    \label{tab:das}
\end{table}

\begin{table}
    \centering
    \begin{tabular}{|c|c|c|c|c|c|c|} 
        \hline
        $N_\mathrm{grid}$  & $T_f$ & $\Delta t$ & Runtime (s) & Runtime (s) & Runtime (s) & Runtime (s) \\ 
        &  &  & (Rosenbrock--Euler) & (EXPRB32) & (EPIRK4s3A) & (EPIRK5P1) \\
        \hline
                                & $2\cdot10^{-7}$ & $1   \cdot \Delta t_\mathrm{CFL}$ & 0.89 & 0.89 & 0.84 & 0.87 \\
        $2^{13} \times 2^{13}$  & $2\cdot10^{-6}$ & $10  \cdot \Delta t_\mathrm{CFL}$ & 0.90 & 0.90 & 0.88 & 0.90 \\ 
                                & $2\cdot10^{-5}$ & $100 \cdot \Delta t_\mathrm{CFL}$ & 0.91 & 0.91 & 0.90 & 0.91 \\
        \hline
                                & $2\cdot10^{-8}$ & $1   \cdot \Delta t_\mathrm{CFL}$ & 0.90 & 0.90 & 0.84 &   \\
        $2^{14} \times 2^{14}$  & $2\cdot10^{-7}$ & $10  \cdot \Delta t_\mathrm{CFL}$ & 0.91 & 0.91 & 0.88 & - \\ 
                                & $2\cdot10^{-6}$ & $100 \cdot \Delta t_\mathrm{CFL}$ & 0.92 & 0.92 & 0.90 &   \\
        \hline
    \end{tabular}
    \caption{Bandwidth achieved on A100, normalised to the maximum practically-achieved bandwidth, for the Burgers' equation with different grid resolutions and time step sizes for selected exponential integrators. All simulations were performed with a user-specified tolerance of $10^{-12}$.}
    \label{tab:bur}
\end{table}

\begin{table}
    \centering
    \begin{tabular}{|c|c|c|c|c|c|} 
        \hline
        Problem & Runtime (s) & Runtime (s) & Runtime (s) & \textbf{Speedup} & \textbf{Speedup} \\
                & 1 CPU core & 112 CPU cores & GPU & 1 CPU core & 112 CPU cores \\
        \hline
        I   &  818.9 &  18.6 & 4.76 & 172.0 & 4.0 \\ 
        \hline
        II  & 1365.4 & 16.6  & 4.14 & 329.8 & 4.0 \\ 
        \hline
        III & 6534.6 & 131.8 & 15.5 & 421.6 & 8.5 \\ 
        \hline
    \end{tabular}
    \caption{Speedups achieved with \texttt{CUDA} over using 1 and 112 core(s) on a CPU for the three problems considered in this work with a configuration of $N_\mathrm{grid} = 2^{14} \times 2^{14}$, $T_f = 2\cdot10^{-6}$, and $\Delta t = 100 \cdot \Delta t_\mathrm{CFL}$). All simulations were performed with a user-specified tolerance of $10^{-12}$.}
    \label{tab:speedups}
\end{table}

Using \lexint, we achieve bandwidths in the range of 1250 - 1325 GB/s for the problems considered in this work. With respect to the theoretical maximum achievable bandwidth of 1555 GB/s on an A1000, we obtain extremely high throughput - on average, between $80\%$ and $85\%$. This efficiency rises to about $90\%$ when compared to the maximum practically-achievable bandwidth. We observe this trend for both linear and nonlinear problems (Tables \ref{tab:das} and \ref{tab:bur}). The speedup achieved by the \texttt{CUDA} implementation compared to the CPU implementation is illustrated in Table \ref{tab:speedups}.

%%%%%%%%%%%%%%%%%%%%%%%%%%%%%%%%%%%%%%%%%%%%%%%%%%%%%%%%%%%%%%%%%%%%%%%%%%%%%%%%%%

\section{Impact and future aspects}
\label{sec:impact}

In our previous work \cite{Deka22_lexint}, we provided an effective implementation of Leja-based exponential integrators in a user-friendly environment. Here, we provide a GPU-based platform for \lexint that enables easy integration into any existing software written in \texttt{CUDA} and \texttt{C++}. \lexint has been highly optimised, both in terms of memory requirements as well as computational speed, and it currently supports simulations on a single GPU. The fundamental structure of the code and the names of the variables are similar in both the \texttt{CUDA} and \texttt{Python} frameworks to enhance readability. 

Simulations of large-scale problems often times use multiple CPUs/GPUs owing to the limiting memory on a single computer system. In the future, we will enable distributing the numerical domain of interest over multiple GPUs and multiple computing nodes by making use of message passing interface. This would be the final step in the full development \lexint that would make it feasible to use it for practical purposes in HPC simulations.

%%%%%%%%%%%%%%%%%%%%%%%%%%%%%%%%%%%%%%%%%%%%%%%%%%%%%%%%%%%%%%%%%%%%%%%%%%%%%%%%%%

\section*{Acknowledgements}
This work is supported, in part, by the Austrian Science Fund (FWF) project id: P32143-N32. This project has received funding, in part, from the European Union's Horizon 2020 research and innovation programme under the Marie Skłodowska-Curie grant agreement No 847476. The views and opinions expressed herein do not necessarily reflect those of the European Commission. The authors acknowledge CINECA and EuroHPC JU for awarding this project access to Leonardo supercomputing hosted at CINECA.

%%%%%%%%%%%%%%%%%%%%%%%%%%%%%%%%%%%%%%%%%%%%%%%%%%%%%%%%%%%%%%%%%%%%%%%%%%%%%%%%%%

% \appendix

% \section{Appendix}
% \label{sec:appendix}

\bibliographystyle{elsarticle-num} 
\bibliography{ref}





% \end{thebibliography}
\end{document}